# Capacity of Control for Stochastic Dynamical Systems Perturbed by Mixed Fractional Brownian Motion with Delay in Control

**Salah H. Abid** *,† **and Uday J. Quaez** †

Department of Mathematics, Education College, Al-Mustansiriya University, Baghdad, Iraq

*,† Correspondence: abidsalah@uomustansiriyah.edu.iq   † uoday1977@gmail.com; Tel.: +964-770-5810-911.

**Abstract:**   In this paper, we discuss the relationships between capacity of control in entropy theory and intrinsic properties in control theory for a class of finite dimensional stochastic dynamical systems described by a linear stochastic differential equations driven by mixed fractional Brownian motion with delay in control. Stochastic dynamical systems can be described as an information channel between the space of control signals and the state space. We study this control to state information capacity of this channel in continuous time. We turned out that, the capacity of control depends on the time of final state in dynamical systems. By using the analysis and representation of fractional Gaussian process, the closed form of continuous optimal control law is derived. The reached optimal control law maximizes the mutual information between control signals and future state over a finite time horizon. The results obtained here are motivated by control to state information capacity for linear systems in both types deterministic and stochastic models that are widely used to understand information flows in wireless network information theory.

The contribution of this paper is that we propose some new relationships between control theory and entropy theoretic properties of stochastic dynamical systems with delay in control. Finally, we present an example that serve to illustrate the relationships between capacity of control and intrinsic properties in control theory.

**Keywords:** controllability; mutual information; capacity of control; entropy; stability; optimal control; stochastic control systems ; delay in control.

## 1. Introduction

One of the fundamental results in the Entropy theory is the capacity of control signals of channels for communication systems. The capacity of control or empowerment which means the maximization of mutual information between the control signals and channel output over all input probability distributions of the channel, plays an important role in many applied fields in engineering, computer science, Economy, physics, chemistry and in other sciences.

Since Shanon in 1948, information theoretic limits of capacity of channels have been studied extensively. The capacities and the capacity of channels of single or multi input signals were studied by Shanon [4], Ahlswede [11] and Liao [6] respectively. Verdu [13] obtained a limiting expression for the capacity regions of multi-access channels with memory. Yu, Rhee and Cioffi [16] characterized the capacity region of a Gaussian multiple access channel with vector inputs and a vector output with or without intersymbol interference. Ranade G. and Sahai A. [12] developed a notation of capacity of control such that gives a fundamental limit (in bits) on the rate at which a controller can dissipate the uncertainty from a system. In a recent important work Tiomkin, Polani and Tishby [14] studied the problem of the control to state information capacity of stochastic dynamic systems driven by Gaussian process with single control in continuous time, when the states are observed only partially and also they derived an efficient solution for computing the control to state information channel. In this paper, we shall study the capacity of control for class of stochastic linear dynamic systems perturbed by mixed fractional Brownian motion with delay in control which are natural generalization of single control. Our work focuses on studying a new relationship between capacity of control in entropy theory and intrinsic properties in control theory for linear systems with delay in control.

The paper is organized as follows: section 2 contains the mathematical model of a linear stochastic systems which described as an information channel between the space of control signals and the state space. Moreover, in this section we define the problem of the information capacity between control signal and observable state with assumption that the system is observable. In section 3, the explicit formula for the mutual information between the resulting state signal at the final time T and the control signals over time t∈ [0, T] is found by using some properties of entropy theory and fractional stochastic analysis theory. We represent and analysis the control process by using orthonormal expansion of Gaussian process in section 4. Section 5 is devoted to study the controllability. In section 6, we derive the optimality conditions for control signal. Relationship between capacity of control and final time is studied in section 7. In section 8, we present a simple numerical example which illustrates the relationship between capacity of control and the intrinsic properties of stochastic dynamic systems.

## 2. System Description

In this paper, we use the following notations: Let $(\Omega, \mathfrak{J}_t, P)$ be a complete probability space with probability measure P on a simple space $\Omega$ defined on a filtration $\{\mathfrak{J}_t, t \in [0,T]\}$. $\mathfrak{J}_t$ is the $\sigma-\text{algebra}$ generated by the random variables $\{W^H_{(s)}, W_{(s)}, s \in [0, t]\}$ and the P-null sets. Let $L^2(\Omega, \mathfrak{J}_t, R^n)$ denote the Hilbert space of all $\mathfrak{J}_t$ measurable square integrable random variables with values in $R^n$.

**Definition (2.1) [15]:**

Let H be a constant belonging to (0,1). A one dimensional fractional Brownian motion $W^H = \{W^H_{(t)}, t \geq 0\}$ of Hurst index H is a continuous and centered Gaussian process with covariance function

$$E(W^H_{(t)} W^H_{(s)}) = \frac{1}{2}(t^{2H} + s^{2H} - |t-s|^{2H}), \text{ for } t, s \geq 0.$$

- If $H = \frac{1}{2}$, then the increments of $W^H$ are non-correlated, and consequently independent. So $W^H$ is a Wiener Process which we denote further by W.
- If $H \in (\frac{1}{2}, 1)$ then the increments are positively correlated.
- If $H \in (0, \frac{1}{2})$ then the increments are negative correlated.

The integral representation of $W^H$ appears as :
$$W^H_{(t)} = \int_0^t K_H(t,s) \, dW_{(s)} \tag{1}$$
where, W is a wiener process and the kernel $K_H(t, s)$ defined as
$$K_H(t,s) = cH \, s^{\frac{1}{2}-H} \int_s^t (u-s)^{H-\frac{3}{2}} u^{H-\frac{1}{2}} \, du \tag{2}$$
$$\frac{\partial K}{\partial t}(t,s) = cH \left(\frac{t}{s}\right)^{H-\frac{1}{2}} (t-s)^{H-\frac{3}{2}} \tag{3}$$

where, $cH = \left[\frac{H(2H-1)}{\beta(2-2H, H-\frac{1}{2})}\right]^{\frac{1}{2}}$, $t > s$ and $\beta$ is a beta function.

In this section, we consider the description of channel for a linear communication system in the following form:
$$\frac{d}{dt}x(t) = Ax(t) + B_1 u(t) + B_2 u(t-h) + G\frac{d}{dt}W(t) + \sigma_1 \frac{d}{dt}W^{H^1}(t) \tag{4}$$
$$y(t) = Dx(t) + \sigma_2 W^{H^2}(t), t \in [0, T] \tag{5}$$
$$x(0) = x_0 \tag{6}$$
$$u(t) = 0, t \in [-h, 0] \tag{7}$$

where, A, G and $\sigma_1$ are constant matrices belongs to $R^{n\times n}$. $B_1$ and $B_2$ are constant matrices belong to $R^{n\times p}$, and the constant matrices $D$ and $\sigma_2$ belong to $R^{p\times n}$. $x_0$ is a centered Gaussian random variable defined on probability space $(\Omega, \Im_0, P)$. The state $x(t) \in C([0,T]; L^2(\Omega, \Im_t, R^n))$ and the values of control $u(t) \in L^2\left(\Omega, \Im_t, L^2(0,T,R^P)\right)$ is a p-dimensional centered Gaussian process defined on probability space $(\Omega, \Im, \{\Im_t\}_{t\geq 0}, P)$ distributed according to probability density function $p(u(t))$, which is restricted to the space of Gaussian process distribution function whose components $u_j$, $j = 1, 2, \ldots, n$ are independent and for any j and final time T the control signal $u_j(t)$, $t \leq T - h$ is independent with $u_j(t)$, $t > T - h$ when h is the time of delay. $\{W_{(t)}, t \in [0,T]\}$ is a multivariate standard Brownian motion (multivariate wiener process) defined on $(\Omega, \Im, \{\Im_t\}_{t\geq 0}, P)$ with values in the space $R^n$ whose components $w_k$, $k = 1, 2, \ldots, n$ are independent. $W^{H^1} = \{W^{H^1}, t \in [0,T]\}$ is a n-multivariate fractional Brownian motion with Hurst parameter $H^1 \in (0,1)^n$ defined in a complete probability space $(\Omega, \Im, \{\Im_t\}_{t\geq 0}, P)$ with values in the space $R^n$ whose components $w_{H^1}^{h_k}$, $k = 1,2,\ldots,n$ are independent. $W^{H^2} = \{W^{H^2}, t \in [0,T]\}$ is a p-multivariate fractional Brownian motion with Hurst parameter $H^2 \in (0,1)^p$ defined in a complete probability space $(\Omega, \Im, \{\Im_t\}_{t\geq 0}, P)$ with values in the space $R^p$ whose components $w_{H^2}^{\tilde{h}_k}$, $k = 1,2,\ldots,n$ are independent.

**Remark (2.2)**

The information channel which is given by the conditional probability distribution $p(y(T)|u(t), t \in [0,T])$ is induced by the stochastic dynamical system $(4-7)$ such that, the system input is perturbed by the noise which resulting by mixed fractional Brownian motion and the system output is perturbed by the noise of fractional Brownian motion.

The linear stochastic dynamic system (4) with initial conditions (6-7) has a unique solution $x(t) \in L^2(\Omega, \Im_t, R^n)$ for any $t \in [0,T]$ see ([2]), which can be represented in the following integral equation:

$$x(t) = e^{At}x_0 + \int_0^t e^{A(t-s)}[B_1 u(s) + B_2 u(s-h)] ds + \int_0^t e^{A(t-s)} G \, dW(s) + \int_0^t e^{A(t-s)} \sigma_1 dW^{H^1}(s) \qquad (8)$$

For any $t \in [-h, 0]$ from the initial control condition (7) the solution of equation (4) with initial condition (5) for $t \in [0, h]$ has the following form:

$$x(t) = e^{At}x_0 + \int_0^t e^{A(t-s)} B_1 u(s) \, ds + \int_0^t e^{A(t-s)} G \, dW(s) + \int_0^t e^{A(t-s)} \sigma_1 dW^{H^1}(t) \qquad (9)$$

Assume that $t > h$, then,

$$x(t) = e^{At}x_0 + \int_0^t e^{A(t-s)} B_1 u(s) \, ds + \int_0^{t-h} e^{A(t-s-h)} B_2 u(s) \, ds + \int_0^t e^{A(t-s)} G \, dW(s) + \int_0^t e^{A(t-s)} \sigma_1 dW^{H^1}(s) \qquad (10)$$

Assume that $T > h$, then the equation (10) is equivalent to the following integral equation:

$$x(t) = e^{At}x_0 + \int_0^{t-h}[e^{A(t-s)} B_1 + e^{A(t-s-h)} B_2]u(s) \, ds + \int_{t-h}^t e^{A(t-s)} B_1 u(s) \, ds \\ + \int_0^t e^{A(t-s)} G \, dW(s) + \int_0^t e^{A(t-s)} \sigma_1 dW^{H^1}(s) \qquad (11)$$

**3. Mutual information**

In this section, we will formulate and obtain the explicit formula for the mutual information between the resulting state signal at the final time T and the control signals over time $t \in [0,T]$, by using some properties of entropy theory and fractional stochastic analysis theory. The mutual information between any two continuous random variables x and y is defined by the difference between the differential entropy of x and the conditional differential entropy of x with respect to y. In other words, the mutual information between x and y defined as $I(x;y) = H(x) - H(x|y)$.

Mathematically, the mutual information between the stochastic processes $y(t)$ and $u(t), t \geq 0$ defined as follows,

$$I(y(T); u(t)) = H(y(T)) - H(y(T)|u(t)), t \in [0,T] \qquad (12)$$

Because the processes $u(t)$, $W(t)$ and $W^{H^1}(t)$ are independent for any $t \in [0,T]$ by assumptions and any one of these process is Gaussian then the equation (11) shows that if the initial condition (3) is multivariate Gaussian random variable or deterministic vector then $x(t)$ is also multivariate Gaussian for any $t \in [0,T]$. The random vector $X \in R^n \sim N(0, S_x)$ with zero mean and covariance matrix S have the differential entropy of X in the following form: $H(X) = \frac{1}{2}\ln(2\pi e)^n \det(S_x)$ (see pp.560 in [1]). From the previous formula, we noted that the differential entropy of multivariate Gaussian random vector depends only on its covariance matrix.

Therefore, the differential entropy of future observed state y(T) at terminal time T is also depending on its covariance matrix $S_y(T)$.

$$H(y(T)) = \frac{1}{2}\ln(2\pi e)^n \det\left(S_y(T)\right) \tag{13}$$

Consequently, under the independence assumption of processes $u(t), W(t), W^{H^1}(t)$ and $W^{H^2}(t)$ then the final state $x(T)$ is independent with multivariate fractional Brownian motion $W^{H^2}$ at time T. Furthermore, from the equation (5) the covariance matrix of y(T) appears as:

$$S_y(T) = DS_x(T)D' + \sigma_2 S_{W^{H^2}}(T)\sigma'_2 \tag{14}$$

Therefore, from equation (11) and the independence assumption of processes $u(t), W(t)$ and $W^{H^1}(t)$ for every $t \in [0, T]$ with initial condition $x_0$ then the covariance matrix of final state $x(T)$ is defined in the following form:

$$S_x(T) = S_{x_0}(T) + S_u(T - h) + S_u(T) + S_{dW}(T) + S_{dW^{H^1}}(T) \tag{15}$$

where,

$S_{x_0}(T) = e^{AT} C_{x_0} e^{A'T}$

$S_u(T - h) = \int_0^{T-h} \int_0^{T-h} [e^{A(T-t_1)}B_1 + e^{A(T-t_1-h)}B_2] \; C_u(t_1, t_2) [e^{A(T-t_2)}B_1 + e^{A(T-t_2-h)}B_2]' dt_1 dt_2$

$S_u(T) = \int_{T-h}^T \int_{T-h}^T [e^{A(T-t_1)}B_1] \; C_u(t_1, t_2) [e^{A(T-t_2)}B_1]' dt_1 dt_2$

$S_{dW}(T) = \int_0^T \int_0^T e^{A(T-t_1)} \; G \; C_{dW}(t_1, t_2) \; G' \; e^{A'(T-t_2)} \; dt_1 dt_2$

$S_{dW^{H^1}}(T) = \int_0^T \int_0^T e^{A(T-t_1)} \sigma_1 C_{dW^{H^1}}(t_1, t_2) \; \sigma_1' \; e^{A'(T-t_2)} \; dt_1 dt_2$

$C_X(t_1, t_2)$ represent the covariance matrix for any two states $x(t_1)$ and $x(t_2)$.

**Lemma (3.1) [9]**

Let V [0,T] be the class of functions such that $f : [0, T] \times \Omega \to R$, f is measurable, $\mathfrak{I}_t$ - adapted and $E\left[\int_0^T (f(t, \omega))^2 dt\right] \leq \infty$. Then for every $f \in V[0, T]$

$$E\left[\int_0^T f(t, \omega) dw(t)\right]^2 = E\left[\int_0^T (f(t, \omega))^2 dt\right] \tag{16}$$

where, w(t) is a standard wiener process.

Consequently, the covariance matrix $S_{dW}(T)$ can be written as

$$S_{dW}(T) = \int_0^T \int_0^T e^{A(T-t_1)} \; G \; G' \; e^{A'(T-t_2)} \; dt_1 dt_2 \tag{17}$$

**Lemma (3.2) [5]**

Let the parameter $\frac{1}{2} < h \leq 1$ then for any functions $\Phi, \varphi \in L^2[0, T] \cap L^1[0, T]$, we have

(i) $\quad E\left(\int_0^T \Phi(t) dw_{(t)}^h \int_0^T \varphi(s) dw_{(s)}^h\right) = h(2h - 1) \int_0^T \int_0^T \Phi(t) \varphi(s) |t - s|^{2h-2} ds dt \tag{18}$

(ii) $\quad E\left(dw_{(t)}^h dw_{(s)}^h\right) = h(2h - 1)|t - s|^{2h-2} ds dt \tag{19}$

Therefore, by using above Lemma and the independence assumption of processes $w_{H^1}^{h_k}(t), k = 1, 2, \dots, n$ for any $t \in [0, T]$, we obtain:

$$S_{dW^{H^1}}(T) = \int_0^T \int_0^T e^{A(T-t_1)} \sigma_1 R(t_1, t_2) \; \sigma_1' e^{A'(T-t_2)} \; dt_1 dt_2 \tag{20}$$

where, $R(t_1, t_2)$ is a n× n diagonal matrix whose kk-th entry is specified by the parameter $h_k$

$r_{kk}(t_1, t_2) = h_k(2h_k - 1)|t_1 - t_2|^{2h_k - 2}$

$r_{kr}(t_1, t_2) = 0, \; r \neq k$

and the covariance function of input signal is defined as

$$C_u(t_1, t_2) = E[u(t_1)u'(t_2)] \tag{21}$$

Now, we will find explicit formula of covariance matrix of $W^{H^2}(T)$. By the independence assumption of the components $w_{H^2}^{\tilde{h}_k}(t), \; k = 1, 2, \dots, p$ for any $t \in [0, T]$, then the covariance matrix $S_{W^{H^2}}(T)$ of $W^{H^2}(T)$ is diagonal, whose kk-th entry is specified by the parameters $\tilde{h}_k$ of the p-multivariate fractional process $W^{H^2}$ at final time T as follows

$$\left[S_{W^{H^2}}(T)\right]_{kk} = T^{2\tilde{h}_k}, \quad k=1,2,\dots,p \tag{22}$$

Consequently,

$$S_y(T) = DS_{x_0}(T)D' + DS_u(T-h)D' + DS_u(T)D' + DS_{dW}(T)D' + DS_{dW^{H^1}}(T)D'$$

$$+\sigma_2 \begin{bmatrix} T^{2\tilde{h}_1} & & 0 \\ & \ddots & \\ 0 & & T^{2\tilde{h}_P} \end{bmatrix}_{P \times P} \sigma'_2 \quad (23)$$

By using the independence assumption of control process $u(t), t \leq T - h$ with the last h of time, we study the mutual information between $y(T)$ and any one of the following control signals $\{u(t), t \leq T - h\}, \{u(t), T - h < t \leq T\}$ and $\{u(t), t \leq T - h\}$.

Therefore, firstly we need to compute the differential entropy of a random vector $y(T)$ conditional on the control process $u(t), t \in [0, T]$.

Assume that $S_{Total}$ is the total covariance of uncontrolled noise, which it comes from the power of mixed fractional Brownian motion. Under the independence assumption of $W(t)$, $W^{H^1}(t)$, $W^{H^2}(t)$ and the initial condition $x_0$, then the total covariance is the following form:

$$S_{Total}(T) = DS_{x_0}(T)D' + DS_{dW}(T)D' + DS_{dW^{H^1}}(T)D' + \sigma_2 \begin{bmatrix} T^{2\tilde{h}_1} & & 0 \\ & \ddots & \\ 0 & & T^{2\tilde{h}_P} \end{bmatrix}_{P \times P} \sigma'_2 \quad (24)$$

The covariance of the future observed state $y(T)$ given the control process $u(t), t \in [0, T]$ is also the total covariance of uncontrolled noise with initial random vector $x_0$

$$S_{Total}(T) = S_{(y|u(t), t \in [0,T])}(T) \quad (25)$$

Consequently,

$$H(y(T)|u(t), t \in [0, T]) = \frac{1}{2}\ln(2\pi e)^n \det(S_{Total}(T)) \quad (26)$$

Therefore, the mutual information between $y(T)$ and the control signals $u(t), t \leq T$ can be written as

$$I(y(T); u(t), t \in [0, T]) = \frac{1}{2}\ln\left\{\det\left(I_{n \times n} + S_{Total}^{-1}(T)(S_u(T-h) + S_u(T))\right)\right\} \quad (27)$$

Now, to find the form of mutual information between $y(T)$ and control signal $u(t), t \leq T - h$ conditioned on the control process $u(t), t \in (T-h, T]$

$$I(y(T)|u(t), t \in (T-h, T]; u(t), t \in [0, T-h]) = H(y(T)|u(t), t \in (T-h, T]) - H(y(T)|u(t), t \in [0, T]) \quad (28)$$

From equation (23) and the independence assumption of $u(t), t \in (T-h, T]$ with $u(t), t \in [0, T-h]$, we have

$S_{y|u(t), t \in (T-h, T]}(T)$

$$= DS_{x_0}(T)D' + DS_u(T-h)D' + DS_{dW}(T)D' + DS_{dW^{H^1}}(T)D' + \sigma_2 \begin{bmatrix} T^{2\tilde{h}_1} & & 0 \\ & \ddots & \\ 0 & & T^{2\tilde{h}_P} \end{bmatrix}_{P \times P} \sigma'_2$$

Hence,

$$H(y(T)|u(t), t \in [T-h, T]) = \frac{1}{2}\ln(2\pi e)^n \det\left(S_{y|u(t), t \in [T-h, T]}(T)\right) \quad (29)$$

Similar, we can to find the mutual information between $y(T)$ and control signal $u(t), t \in (T-h, T]$ conditioned on the control process $u(t), t \in [0, T-h]$ from the following equation

$S_{y|u(t), t \in [0, T-h]}(T)$

$$= DS_{x_0}(T)D' + DS_u(T)D' + DS_{dW}(T)D' + DS_{dW^{H^1}}(T)D' + \sigma_2 \begin{bmatrix} T^{2\tilde{h}_1} & & 0 \\ & \ddots & \\ 0 & & T^{2\tilde{h}_P} \end{bmatrix}_{P \times P} \sigma'_2 \quad (30)$$

Where,

$$I(y(T)|u(t), t \in [0, T-h]; u(t), t \in (T-h, T]) = H(y(T)|u(t), t \in [0, T]) - H(y(T)|u(t), t \in (T-h, T]) \quad (31)$$

### 4. Representation and Analysis of Control Process

In this section, we represent and analysis the control process by using orthonormal expansion of Gaussian process.

**Lemma (4.1) (Karhunen-Loeve Theorem ) [15]**

Let $x(t), t \in [a, b]$ be a central square integrable stochastic process defined over a probability space $(\Omega, \Im, P)$ with continuous covariance function $C_x$ then $x(t) = \sum_{j=1}^{\infty} z_j e_j(t)$, where, $e_j, j = 1,2,...$ is a set of orthonormal basis on $L^2([a, b])$ and the random variables $z_j, j = 1,2,...$ are independent and have zero mean. From this above lemma for any component $u_k(t)$, k= 1,2,…,p of the control signal $u(t) \in L^2(\Omega, \Im_t, L^2(0, T - h, R^P))$ may be represented by appropriate choice of $\{e_{j,k}(t)\}_{j=1}^{\infty}$, such that the set of functions $\{e_{j,k}(t)\}_{j=1}^{\infty}$ is countable of real orthonormal functions in $L^2([0, T - h])$ and $u_{j,k}, j = 1,2,...$ is a sequence of independent Gaussian random variables as follows

$$u_k(t) = \sum_{j=1}^{\infty} u_{j,k} e_{j,k}(t) \qquad (32)$$

Also, for any component $u_k(t)$, k= 1,2,…,p of the control signal $u(t) \in L^2(\Omega, \Im_t, L^2(T - h, T, R^P))$ can be written as

$$u_k(t) = \sum_{j=1}^{\infty} \hat{u}_{j,k} \theta_{j,k}(t) \qquad (33)$$

Such that $\{\theta_{j,k}(t)\}_{j=1}^{\infty}$ is a countable of real orthonormal functions in $L^2(T - h, T, R^P)$ and $\hat{u}_{j,k}, j = 1,2,...$ is a sequence of independent Gaussian random variables.

Now, for any $t \in [0, T - h]$, the p-dimensional control process can be represented as follows:

$$u(t) = \begin{bmatrix} \sum_{j=1}^{\infty} u_{j,1} e_{j,1}(t) \\ \sum_{j=1}^{\infty} u_{j,2} e_{j,2}(t) \\ \cdot \\ \cdot \\ \cdot \\ \sum_{j=1}^{\infty} u_{j,p} e_{j,p}(t) \end{bmatrix} \qquad (34)$$

where, $\{e_{j,k}(t)\}_{j=1}^{\infty}$ is countable of real orthonormal functions in $L^2(0, T - h, R^P)$.

By independence assumption of components of control process $u(t)$ then its covariance function is diagonal, whose kk-th entry is specified by the parameters $\sigma_{jk}$ of the control process as follows

$$[C_u(t_1, t_2)]_{kk} = \sum_{j=1}^{\infty} \sigma_{jk} e_{j,k}(t_1) e_{j,k}(t_2), k=1,2,\ldots \qquad (35)$$

Also, for any $t \in [T - h, T]$ the p-dimensional control process can be represented as:

$$u(t) = \begin{bmatrix} \sum_{j=1}^{\infty} \hat{u}_{j,1} \theta_{j,1}(t) \\ \sum_{j=1}^{\infty} \hat{u}_{j,2} \theta_{j,2}(t) \\ \cdot \\ \cdot \\ \cdot \\ \sum_{j=1}^{\infty} \hat{u}_{j,p} \theta_{j,p}(t) \end{bmatrix} \qquad (36)$$

and

$$[C_u(t_1, t_2)]_{kk} = \sum_{j=1}^{\infty} \omega_{jk} \theta_{j,k}(t_1) \theta_{j,k}(t_2), k=1,2,\ldots \qquad (37)$$

Suppose that $B = B_1 + e^{-Ah} B_2$ then the covariance of $u(t), t \leq T - h$ appearers as

$$S_u(T - h) = \int_0^{T-h} \int_0^{T-h} e^{A(T-t_1)} \sum_{k=1}^{p} \sum_{j=1}^{\infty} b_k \sigma_{jk} e_{j,k}(t_1) e_{j,k}(t_2) b_k' e^{A'(T-t_2)} dt_1 dt_2 \qquad (38)$$

where, $b_k = B[:, k]$ is a k-th column of a matrix $B$. Similarity,

$$S_u(T) = \int_{T-h}^{T} \int_{T-h}^{T} e^{A(T-t_1)} \sum_{k=1}^{p} \sum_{j=1}^{\infty} b_{1k} \omega_{jk} \theta_{j,k}(t_1) \theta_{j,k}(t_2) b_{1k}' e^{A'(T-t_2)} dt_1 dt_2 \qquad (39)$$

Where, $b_{1k} = B_1[:, k]$ is a k-th column of matrix $B_1$

By using the assumption that $u(t) \in L^2(\Omega, \Im_t, L^2(0, T, R^P))$, then

$$\int_0^{T-h} \|u(t)\|^2_{L^2(\Omega,\Im_t,L^2(0,T,R^P))} dt \leq M_1 \qquad (40)$$

$$\int_{T-h}^{T} \|u(t)\|^2_{L^2(\Omega,\Im_t,L^2(0,T,R^P))} dt \leq M_2 \qquad (41)$$

$M_1$ and $M_2$ are constants.

Consequently,

$$\text{Tr}\left\{\int_0^{T-h} C_u(t, t) dt\right\} = \sum_{k=1}^{p} \sum_{j=1}^{\infty} \int_0^{T-h} \sigma_{jk} e_{j,k}(t) e_{j,k}(t) dt = \sum_{k=1}^{p} \sum_{j=1}^{\infty} \sigma_{jk} \leq M_1 \qquad (42)$$

$$\text{Tr}\left\{\int_{T-h}^{T} C_u(t, t) dt\right\} = \sum_{k=1}^{p} \sum_{j=1}^{\infty} \int_{T-h}^{T} \omega_{jk} \theta_{j,k}(t) \theta_{j,k}(t) dt = \sum_{k=1}^{p} \sum_{j=1}^{\infty} \omega_{jk} \leq M_2 \qquad (43)$$

## 5. Controllability

In this section we study the controllability of the system (4-7) by using the properties of deterministic linear system with delay in control and controllable matrix. Consider the deterministic linear system with delay control as:

$$\frac{d}{dt}x(t) = Ax(t) + B_1 u(t) + B_2 u(t-h), t \in [0, T] \tag{44}$$

$$x(0) = x_0 \tag{45}$$

$$u(t) = 0, t \in [-h, 0] \tag{46}$$

Define the linear control operator $L_0^T u : L^2\big(\Omega, \mathfrak{I}_t, L^2(0, T, R^P)\big) \to L^2\big(\Omega, \mathfrak{I}_t, L^2(0, T, R^n)\big)$ as

$$L_0^T u = \int_0^{T-h} [e^{A(T-s)} B_1 + e^{A(T-s-h)} B_2] u(s)\, ds + \int_{T-h}^T e^{A(T-s)} B_1 u(s)\, ds \tag{47}$$

It is clear that, $L_0^T u$ is a bounded. The adjoint operator

$(L_0^T u)^* : L^2\big(\Omega, \mathfrak{I}_t, L^2(0, T, R^n)\big) \to L^2\big(\Omega, \mathfrak{I}_t, L^2(0, T, R^P)\big)$ is defined by

$$(L_0^T)^* z = \begin{cases} \big(e^{A(T-t)} B\big)' E(z) & , t \in [0, T-h] \\ \big(e^{A(t-t)} B_1\big)' E(z) & , t \in [T-h, T] \end{cases} \tag{48}$$

Also, we defined the controllability operator $G_T$ associated with control operator $L_0^T$ in deterministic case of equation (44) as

$$G_T = L_0^T (L_0^T)^* = \int_0^{T-h} e^{A(T-t)} BB' e^{A'(T-t)} dt + \int_{T-h}^T e^{A(T-t)} B_1 B_1' e^{A'(T-t)} dt \tag{49}$$

**Lemma (5.1) [7]**

The following conditions are equivalent

*i* The deterministic control system (44-46) is controllable on $[0, T]$.

*ii* The controllability matrix $G_T$ in (49) is nonsingular.

*iii* The matrix $[B_1, B_2, AB_1, AB_2, \ldots, A^{n-1} B_1, A^{n-1} B_2]$ is full rank.

**Lemma (5.2)**

The following conditions are equivalent

a. The deterministic control system (44-46) is controllable on $[0, T]$.

b. The stochastic control system (4-7) is controllable on $[0, T]$.

Proof: similar of the proof of theorem (1) in [7].

Note that, the following notations will be used in this paper

$$g_{A,b_k}(t) = e^{A(T-t)} b_k \in R^n \tag{50}$$

$$g_{A,b_{1k}}(t) = e^{A(T-t)} b_{1k} \in R^n \tag{51}$$

$$v_{ik}(T) = \int_0^{T-h} \big(g_{A,b_k}(t) e_{j,k}(t)\big) dt \tag{52}$$

$$z_{ik}(T) = \int_{T-h}^T g_{A,b_{1k}}(t) \theta_{j,k}(t) dt \tag{53}$$

Consequently, the covariance matrices $S_u(T)$ and $S_u(T-h)$ can be written as

$$S_u(T) = \sum_{k=1}^p \sum_{j=1}^\infty \omega_{jk}\, z_{ik}(T) z'_{ik}(T) \tag{54}$$

$$S_u(T-h) = \sum_{k=1}^p \sum_{j=1}^\infty \sigma_{jk}\, v_{ik}(T) v'_{ik}(T) \tag{55}$$

$$G_T(b_k, b_{1k}) = \int_0^{T-h} g_{A,b_k}(t) g'_{A,b_k}(t) dt + \int_{T-h}^T g_{A,b_{1k}}(t) g'_{A,b_{1k}}(t) dt \in R^{n \times n} \tag{56}$$

Assume that the deterministic linear control system which corresponding to stochastic control system (4-7) is controllable. This mean that the matrix $G_T(b_k, b_{1^k})$ is positive definite for any k=1,2,...,p.

In the next section, we shall formulate and prove the conditions for optimal control. In other words, we find the expression of control signal which maximizes the mutual information between control signals and future observe stat for the stochastic dynamical system (4-7).

## 6. Optimality

Consider the following constraint optimization :

$$\max_{p(u(t), t \in [0,T])} I(y(T); u(t), t \in [0, T])$$

Subject to

$$\frac{d}{dt} x(t) = Ax(t) + B_1 u(t) + B_2 u(t-h) + G \frac{d}{dt} W(t) + \sigma_1 \frac{d}{dt} W^{H^1}(t)$$

$$y(t) = Dx(t) + \sigma_2 W^{H^2}(t) \quad , t \in [0, T]$$

$$x(0) = x_0$$

$$u(t) = 0 \quad , t \in [-h, 0] \tag{57}$$

By using the representation of control process in section 4, the equation (27) can be written in the following form:

$$I(y(T); u(t), t \in [0, T])$$
$$= \frac{1}{2} \ln\{\det(I_{n \times n} + S_{Total}^{-1}(T) \left( \sum_{k=1}^{p} \sum_{j=1}^{\infty} \sigma_{jk} D v_{ik}(T) v'_{ik}(T) D' + \sum_{k=1}^{p} \sum_{j=1}^{\infty} \omega_{jk} D z_{ik}(T) z'_{ik}(T) D' \right) \tag{58}$$

Constraint conditions (42-43) with orthonormality conditions implies to a new expression for optimization problem (57)

$$\max_{\sigma, \omega, \{\theta_{j,k}(t)\}_{jk}, \{e_{j,k}(t)\}_{jk}} \frac{1}{2} \ln\{\det(I_{n \times n} + S_{Total}^{-1}(T)$$
$$\cdot \left( \sum_{k=1}^{p} \sum_{j=1}^{\infty} \sigma_{jk} D v_{ik}(T) v'_{ik}(T) D' + \sum_{k=1}^{p} \sum_{j=1}^{\infty} \omega_{jk} D z_{ik}(T) z'_{ik}(T) D' \right) \tag{59}$$

Subject to

$$\sum_{k=1}^{p} \sum_{j=1}^{\infty} \sigma_{jk} = M_1 \tag{60}$$

$$\sum_{k=1}^{p} \sum_{j=1}^{\infty} \omega_{jk} = M_2 \tag{61}$$

$$\sigma_{jk} \geq 0 \tag{62}$$

$$\omega_{jk} \geq 0 \tag{63}$$

$$\int_0^{T-h} e_{j,k}(t) \, e_{i,k}(t) dt = \delta_{ji} \tag{64}$$

$$\int_{T-h}^{T} \theta_{j,k}(t) \, \theta_{i,k}(t) dt = \delta_{ji} \tag{65}$$

where,

$$\delta_{ji} = \begin{cases} 0 & , \text{if } i = j \\ 1 & , \text{if } i \neq j \end{cases} \quad , \sigma = \{\sigma_{jk}\}_{jk} \text{ and } \omega = \{\omega_{jk}\}_{jk}$$

Therefore, the Lagrange function can be written as

$$L = \frac{1}{2} \ln\{\det(I_{n \times n} + S_{Total}^{-1}(T) \cdot \left( \sum_{k=1}^{p} \sum_{j=1}^{\infty} \sigma_{jk} D v_{ik}(T) v'_{ik}(T) D' + \sum_{k=1}^{p} \sum_{j=1}^{\infty} \omega_{jk} D z_{ik}(T) z'_{ik}(T) D' \right)$$
$$- \gamma_1 \left( \sum_{k=1}^{p} \sum_{j=1}^{\infty} \sigma_{jk} - M_1 \right) - \gamma_2 \left( \sum_{k=1}^{p} \sum_{j=1}^{\infty} \omega_{jk} - M_2 \right) - \sum_{k=1}^{p} \sum_{j=1}^{\infty} \lambda_{jk} \sigma_{jk} - \sum_{k=1}^{p} \sum_{j=1}^{\infty} \mu_{jk} \omega_{jk}$$
$$- \sum_{k=1}^{p} \sum_{j,i=1}^{\infty} \beta_{k,j,i} \left( \int_0^{T-h} e_{j,k}(t) e_{i,k}(t) dt - \delta_{ji} \right) - \sum_{k=1}^{p} \sum_{j,i=1}^{\infty} \alpha_{k,j,i} \left( \int_{T-h}^{T} \theta_{j,k}(t) \theta_{i,k}(t) dt - \delta_{ji} \right) \tag{66}$$

where, $\gamma_1, \gamma_2, \lambda_{jk}, \mu_{jk}, \beta_{k,j,i}$ and $\alpha_{k,j,i}$ are the Lagrange multipliers.

Consequently, the corresponding KKT (Karush, Kuhn and Tucker [8]) optimality conditions of optimization problem (59-65) are in the following form, for any $t \in [0, T]$, k=1,2,...,p and j, i=1,2,...

(i) $\frac{\delta L}{\delta e_{j,k}}(t) = 0$ (67)

$\frac{\delta L}{\delta \theta_{j,k}}(t) = 0$ (68)

$\frac{\partial L}{\partial \sigma_{jk}} = 0$ (69)

$\frac{\partial L}{\partial \omega_{jk}} = 0$ (70)

(ii) $\sum_{k=1}^{p} \sum_{j=1}^{\infty} \sigma_{jk} - M_1 = 0$ (71)

$\sum_{k=1}^{p} \sum_{j=1}^{\infty} \omega_{jk} - M_2 = 0$ (72)

$\sigma_{jk} \geq 0$ (73)

$\omega_{jk} \geq 0$ (74)

$\int_0^T e_{j,k}(t) e_{i,k}(t) dt - \delta_{ji} = 0$ (75)

$\int_{T-h}^T \theta_{j,k}(t) \theta_{i,k}(t) dt - \delta_{ji} = 0$ (76)

(iii) $\lambda_{jk} \sigma_{jk} = 0$ (77)

$\mu_{jk} \omega_{jk} = 0$ (78)

(iv) $\lambda_{jk} \geq 0, \ \mu_{jk} \geq 0$ (79)

Clearly, the nonlinear optimization problem (59-65) is a convex with respect to $\{\sigma_{jk}\}_{jk}$ and $\{\omega_{jk}\}_{jk}$ for a given set of the expansion functions $\{e_{j,k}(t)\}_{jk}$ and $\{\theta_{j,k}(t)\}_{jk}$. The partial control signal without the rm-th control component for Gaussian process representation in (32-33) is defined as

$u_{\widetilde{rm}}(t) = \sum_{jk \neq rm} u_{j,k} e_{j,k}(t) \quad , t \in [0, T-h]$ (80)

$u_{\widetilde{rm}}(t) = \sum_{jk \neq rm} \hat{u}_{j,k}(t)_{j,k} \theta_{j,k}(t) \quad , t \in (T-h, T]$ (81)

Therefore, the covariance matrices of control signals (38-39) can be written as:

$S_u(T-h) = \sum_{jk \neq rm} \sigma_{jk} v_{jk}(T) v'_{jk}(T) + \sigma_{rm} v_{rm}(T) v'_{rm}(T)$ (82)

$S_u(T) = \sum_{jk \neq rm} \omega_{jk} z_{jk}(T) z'_{jk}(T) + \omega_{rm} z_{rm}(T) z'_{rm}(T)$ (83)

Consequently, the covariance matrix of partial control signal $u_{\widetilde{rm}}(t)$ appears as:

$S_{u_{\widetilde{rm}}}(T-h) = S_u(T-h) - \sigma_{rm} v_{rm}(T) v'_{rm}(T)$ (84)

$S_{u_{\widetilde{rm}}}(T) = S_u(T) - \omega_{rm} z_{rm}(T) z'_{rm}(T)$ (85)

The covariance matrix of a final observable state conditional by rm-th control process component can be written as:

$S_{(y|u_{rm})}(T) = D S_{u_{\widetilde{rm}}}(T-h) D' + D S_{u_{\widetilde{rm}}}(T) D' + S_{Total}(T)$ (86)

Substituting $S_{Total}(T)$ in Lagrange function (81), we get

$L = \frac{1}{2} \ln \det \{S_{(y|u_{rm})}(T) - D S_{u_{\widetilde{rm}}}(T-h) D' - D S_{u_{\widetilde{rm}}}(T) D' \sum_{k=1}^{p} \sum_{j=1}^{\infty} \sigma_{jk} D v_{ik}(T) v'_{ik}(T) D'$
$+ \sum_{k=1}^{p} \sum_{j=1}^{\infty} \omega_{jk} D z_{ik}(T) z'_{ik}(T) D'\} - \frac{1}{2} \ln \det (S_{Total}(T)) - \gamma_1 (\sum_{k=1}^{p} \sum_{j=1}^{\infty} \sigma_{jk} - M_1)$
$- \gamma_2 (\sum_{k=1}^{p} \sum_{j=1}^{\infty} \omega_{jk} - M_2) - \sum_{k=1}^{p} \sum_{j=1}^{\infty} \lambda_{jk} \sigma_{jk} - \sum_{k=1}^{p} \sum_{j=1}^{\infty} \mu_{jk} \omega_{jk}$
$- \sum_{k=1}^{p} \sum_{j,i=1}^{\infty} \beta_{k,j,i} \left( \int_0^{T-h} e_{j,k}(t) e_{i,k}(t) dt - \delta_{ji} \right) - \sum_{k=1}^{p} \sum_{j,i=1}^{\infty} \alpha_{k,j,i} \left( \int_{T-h}^{T} \theta_{j,k}(t) \theta_{i,k}(t) dt - \delta_{ji} \right)$ (87)

Using equations (84-85), we get

$L = \{\frac{1}{2} \ln \det \{S_{(y|u_{rm})}(T) + \sigma_{rm} D v_{rm}(T) v'_{rm}(T) D' + \omega_{rm} D z_{rm}(T) z'_{rm}(T) D'\}$
$- \frac{1}{2} \ln \det (S_{Total}(T)) - \gamma_1 (\sum_{k=1}^{p} \sum_{j=1}^{\infty} \sigma_{jk} - M_1) - \gamma_2 (\sum_{k=1}^{p} \sum_{j=1}^{\infty} \omega_{jk} - M_2)$
$- \sum_{k=1}^{p} \sum_{j=1}^{\infty} \lambda_{jk} \sigma_{jk} - \sum_{k=1}^{p} \sum_{j=1}^{\infty} \mu_{jk} \omega_{jk} - \sum_{k=1}^{p} \sum_{j,i=1}^{\infty} \beta_{k,j,i} \left( \int_0^{T-h} e_{j,k}(t) e_{i,k}(t) dt - \delta_{ji} \right)$
$- \sum_{k=1}^{p} \sum_{j,i=1}^{\infty} \alpha_{k,j,i} \left( \int_{T-h}^{T} \theta_{j,k}(t) \theta_{i,k}(t) dt - \delta_{ji} \right)$ (88)

Now, to compute the ordinary derivative of Lagrange function with respect to $\sigma_{rm}$ for each r =1,2,... and m =1,2,...,p

$\frac{\partial L}{\partial \sigma_{rm}} = \frac{v'_{rm}(T) D' S^{-1}_{(y|u_{rm})}(T) D v_{rm}(T)}{1 + \sigma_{rm} v'_{rm}(T) D' S^{-1}_{(y|u_{rm})}(T) D v_{rm}(T) + \omega_{rm} z'_{rm}(T) D' S^{-1}_{(y|u_{rm})} D z_{rm}(T)} - \gamma_1 - \lambda_{rm}$ (89)

Equating the ordinary derivative of Lagrange function to zero for each r =1, 2,... , m =1,2,...,p, and using the KKT conditions (77-79), we have $\sigma_{rm} = 0$ leads to

$\frac{v'_{rm}(T) D' S^{-1}_{(y|u_{rm})}(T) D v_{rm}(T)}{1 + \omega_{rm} z'_{rm}(T) D' S^{-1}_{(y|u_{rm})} D z_{rm}(T)} - \gamma_1 - \lambda_{rm} = 0$ (90)

and $\sigma_{rm} > 0$ leads to

$$\frac{v'_{rm}(T)D'\,S^{-1}_{(y|u_{rm})}(T)Dv_{rm}(T)}{1+\sigma_{rm}Dv_{rm}(T)S^{-1}_{(y|u_{rm})}v'_{rm}(T)D'+\omega_{rm}z'_{rm}(T)D'S^{-1}_{(y|u_{rm})}Dz_{rm}(T)} - \gamma_1 = 0 \quad (91)$$

Hence, for any k=1,2,…,p and j, i=1,2,… , we get

$$\sigma_{jk} = max\left\{0, \frac{1}{\gamma_1} - \frac{1}{v'_{jk}(T)D'\,S^{-1}_{(y|u_{jk})}(T)Dv_{jk}(T)} - \frac{\omega_{jk}z'_{jk}(T)D'S^{-1}_{(y|u_{jk})}(T)\,Dz_{jk}(T)}{v'_{jk}(T)D'\,S^{-1}_{(y|u_{jk})}(T)Dv_{jk}(T)}\right\} \quad (92)$$

From KKT condition (71), we have

$$\sum_{k=1}^{p}\sum_{j=1}^{\infty} max\left\{0, \frac{1}{\gamma_1} - \frac{1}{v'_{jk}(T)D'\,S^{-1}_{(y|u_{jk})}(T)Dv_{jk}(T)} - \frac{\omega_{jk}z'_{jk}(T)D'S^{-1}_{(y|u_{jk})}(T)\,Dz_{jk}(T)}{v'_{jk}(T)D'\,S^{-1}_{(y|u_{jk})}(T)Dv_{jk}(T)}\right\} = M_1 \quad (93)$$

Similarity, the ordinary derivative of Lagrange function with respect to $\omega_{rm}$ for each r =1,2,… and m =1,2,…,p is,

$$\frac{\partial L}{\partial \omega_{rm}} = \frac{z'_{rm}(T)D'\,S^{-1}_{(y|u_{rm})}(T)Dz_{rm}(T)}{1+\sigma_{rm}v'_{rm}(T)D'S^{-1}_{(y|u_{rm})}(T)Dv_{rm}(T)+\omega_{rm}z'_{rm}(T)D'S^{-1}_{(y|u_{rm})}Dz_{rm}(T)} - \gamma_2 - \mu_{rm} \quad (94)$$

$\omega_{rm} = 0$ implies

$$\frac{z'_{rm}(T)D'\,S^{-1}_{(y|u_{rm})}(T)Dz_{rm}(T)}{1+\sigma_{rm}v'_{rm}(T)D'S^{-1}_{(y|u_{rm})}(T)Dv_{rm}(T)+} - \gamma_2 - \mu_{rm} = 0 \quad (95)$$

$\omega_{rm} > 0$ implies

$$\frac{z'_{rm}(T)D'\,S^{-1}_{(y|u_{rm})}(T)Dz_{rm}(T)}{1+\sigma_{rm}Dv_{rm}(T)S^{-1}_{(y|u_{rm})}v'_{rm}(T)D'+\omega_{rm}z'_{rm}(T)D'S^{-1}_{(y|u_{rm})}Dz_{rm}(T)} - \gamma_2 = 0 \quad (96)$$

Therefore, for any k=1,2,…,p and j, i=1,2,… , we get

$$\omega_{jk} = max\left\{0, \frac{1}{\gamma_2} - \frac{1}{z'_{jk}(T)D'\,S^{-1}_{(y|u_{jk})}(T)Dz_{jk}(T)} - \frac{\sigma_{jk}Dv_{jk}(T)S^{-1}_{(y|u_{jk})}v'_{jk}(T)D'}{z'_{jk}(T)D'\,S^{-1}_{(y|u_{jk})}(T)Dz_{jk}(T)}\right\} \quad (97)$$

and

$$\sum_{k=1}^{p}\sum_{j=1}^{\infty} max\left\{0, \frac{1}{\gamma_2} - \frac{1}{z'_{jk}(T)D'\,S^{-1}_{(y|u_{jk})}(T)Dz_{jk}(T)}\right\} - \frac{\sigma_{jk}Dv_{jk}(T)S^{-1}_{(y|u_{jk})}v'_{jk}(T)D'}{z'_{jk}(T)D'\,S^{-1}_{(y|u_{jk})}(T)Dz_{jk}(T)} = M_2 \quad (98)$$

Now, we derive the optimality conditions for the expansion functions $\{e_{j,k}(t)\}$ by computing the derivative of the Lagrange function (66). A variation of any function g can be represented as $\delta g(t) = \varepsilon \varphi(t)$, where the quantity $\varepsilon$ is an infinitesimal number and $\varphi$ is a test function. The ordinary or partial derivative of functional G(g) can be defined via the variation $\delta G$, which results from variation of g by $\delta G$

$$\delta G = G(g + \delta g) - G(g) \quad (99)$$

Using the Taylor expansion of $G(g + \delta g) = G(g + \varepsilon \varphi)$ we get

$$G(g + \varepsilon\varphi) = G(g) + \frac{d}{d\varepsilon}G(g + \varepsilon\varphi)\Big|_{\varepsilon=0}\varepsilon + \frac{1}{2}\frac{d^2}{d\varepsilon^2}G(g + \varepsilon\varphi)\Big|_{\varepsilon=0}\varepsilon^2 + \cdots \quad (100)$$

The functional derivative G with respect to g is defined as:

$$\frac{d}{dg}G(g + \varepsilon\varphi)\Big|_{\varepsilon=0} = \int \frac{\delta G(g)}{\delta g(t)}\varphi(t)dt \quad (101)$$

Consequently,

$$\frac{\delta L}{\delta e_{r,m}} = Tr\left\{S_y^{-1}(T)\sigma_{rm}D\left(g_{A,b_m}(t)v'_{rm}(T) + v_{rm}(T)g'_{A,b_m}(t)\right)D'\right\} + 2\sum_{i=1}^{\infty}\beta_{m,r,i}e_{i,m}(t) \quad (102)$$

The optimality condition (67) implies to the following equation:

$$\sigma_{rm}Tr\{S_y^{-1}(T)Dg_{A,b_m}(t)v'_{rm}(T)\,D'\} + \sigma_{rm}Tr\{S_y^{-1}(T)Dv_{rm}(T)g'_{A,b_m}(t)\,D'\} + 2\sum_{i=1}^{\infty}\beta_{m,r,i}e_{i,m}(t) = 0 \quad (103)$$

Using the property that the trace of matrices is symmetric and the symmetry of $S_y^{-1}(T)$, we obtain

$$\sigma_{rm}Tr\{S_y^{-1}(T)Dv_{rm}(T)g'_{A,b_m}(t)\,D'\} = -\sum_{i=1}^{\infty}\beta_{m,r,i}e_{i,m}(t) \quad (104)$$

This equation is true for any k =1, 2,…, p, j =1,2,… and t ∈ [0, T − h]

$$\sigma_{jk}Tr\{S_y^{-1}(T)Dv_{jk}(T)g'_{A,b_k}(t)\,D'\} = -\sum_{i=1}^{\infty}\beta_{k,j,i}e_{i,k}(t) \quad (105)$$

Multiplying both sides by $e_{\tau,k}(t)$, τ=1, 2,… and taking an integral over t ∈ [0, T − h], we have

$$\sigma_{jk} \, \text{Tr}\left(S_y^{-1}(T) D v_{jk}(T) \int_0^{T-h} e_{\tau,k}(t) g'_{A,b_k}(t) dt \, D'\right) = -\sum_{i=1}^{\infty} \beta_{k,j,i} \int_0^{T-h} e_{\tau,k}(t) \, e_{i,k}(t) dt \tag{106}$$

Using orthonormality of a set of the functions $\{e_{j,k}(t)\}$, we get

$$\sigma_{jk} \, \text{Tr}(S_y^{-1}(T) D v_{jk}(T) v'_{\tau k}(T) \, D') = -\beta_{k,j,\tau} \tag{107}$$

Substituting $\beta_{k,j,\tau}$ in equation (106), we get

$$\sigma_{jk} \text{Tr}\{S_y^{-1}(T) D v_{jk}(T) g'_{A,b_k}(t) \, D'\} = \sum_{i=1}^{\infty} \sigma_{jk} \, \text{Tr}\left(S_y^{-1}(T) D v_{jk}(T) v'_{ik}(T) \, D'\right) e_{i,k}(t) \tag{108}$$

Multiplying above equation by $g_{A,b_k}(t)$ and taking an integral over $t \in [0, T-h]$, we have

$$\sigma_{jk} \int_0^{T-h} g_{A,b_k}(t) \, g'_{A,b_k}(t) dt \, D' S_y^{-1}(T) \, D v_{jk}(T) = \sigma_{jk} \sum_{i=1}^{\infty} v_{ik}(T) \, v'_{ik}(T) \, D' S_y^{-1}(T) D v_{jk}(T) \tag{109}$$

Therefore, equation (109) can be rewritten as

$$\sigma_{jk} \left[\int_0^{T-h} g_{A,b_k}(t) \, g'_{A,b_k}(t) dt - \sum_{i=1}^{\infty} v_{ik}(T) \, v'_{ik}(T)\right] . D' S_y^{-1}(T) \, D v_{jk}(T) = 0 \tag{110}$$

Assume that

$$\widehat{G}_1 = G_1(T) - \sum_{i=1}^{\infty} v_{ik}(T) \, v'_{ik}(T) \tag{111}$$

where, $G_1(b_k) = \int_0^{T-h} g_{A,b_k}(t) \, g'_{A,b_k}(t) dt$

The equation (110) is satisfied if at least one of the following cases is hold

i. $\sigma_{jk} = 0$
ii. $v_{jk}(T) = 0$
iii. $\widehat{G}_1 = 0$
iv. $D' S_y^{-1}(T) \, D v_{jk}(T)$ belong to $\text{null}(\widehat{G}_1)$.

**Lemma (6.1)**

Assume that the matrix D is nonsingular and there exist r, m such that $\sigma_{rm} > 0$, but $v_{rm}(T) = 0$ then $I(y(T)|u(t), t \in (T-h,T]; u(t), t \in [0, T-h]) = 0$

**Proof:**

Suppose that there exist r, m such that $\sigma_{rm} > 0$ and $v_{rm}(T) = 0$. Since $S_{(y|u_{rm})}^{-1}(T)$ is a positive definite and $v_{rm}(T) = 0$ then

$$v'_{rm}(T) D' S_{(y|u_{rm})}^{-1}(T) \, D v_{rm}(T) = 0 \tag{112}$$

But $\sigma_{rm} > 0$ by assumption of this lemma, then from equation (91), we get $\gamma_1 = 0$.

Consequently, for all j, k

$$v'_{jk}(T) D' \, S_{(y|u_{jk})}^{-1}(T) D v_{jk}(T) = 0 \tag{113}$$

This means that $v_{jk}(T) = 0$, in other words,

$$\sum_{j,k:\sigma_{jk}>0} \sigma_{jk} \, D v_{jk}(T) v'_{jk}(T) D' = 0 \tag{114}$$

If $\sigma_{jk} = 0$ then

$$\sum_{j,k:\sigma_{jk}=0} \sigma_{jk} \, D v_{jk}(T) v'_{jk}(T) D' = 0 \tag{115}$$

Adding equation (115) to equation (114), we get

$$\sum_{j,k} \sigma_{jk} \, D v_{jk}(T) v'_{jk}(T) D' = 0 \tag{116}$$

This completes the proof.

Now, since there exists at least one of a set $\{\sigma_{jk}\}_{jk}$ which is greater than zero. Therefore, if $\sigma_{rm} > 0$ for some r, m, then from the above lemma, the mutual information between y(T) conditional by $u(t), t \in (T-h,T]$ and the control signal $u(t), t \in [0, T-h]$ becomes zero and $v_{jk}(T) = 0$, for any j, k when $v_{rm}(T) = 0$. Therefore, $v_{rm}(T) \neq 0$. Since $S_y^{-1}(T)$ is a positive definite, then $D' S_y^{-1}(T) \, D v_{jk}(T) \neq 0$. In the third case, the solution of equation (110) is a decomposition of the matrix $G_1(b_k)$ in to the sum of one rank matrices

$$G_1(b_k) = \sum_{i=1}^{\infty} v_{ik}(T) \, v'_{ik}(T) \tag{117}$$

Similarity,

$$G_2(b_{1k}) = \sum_{i=1}^{\infty} z_{ik}(T) \, z'_{ik}(T) \tag{118}$$

where,
$$G_2(b_{1k}) = \int_{T-h}^{T} g_{A,b_{1k}}(t) \, g'_{A,b_{1k}}(t) \, dt \tag{119}$$

**Theorem (6.1)**

Assume that both matrices $G_1(b_k)$ and $G_2(b_{1k})$ are positive definite then the optimal sets $\{\sigma_{jk}\}_{jk}$ and $\{\omega_{jk}\}_{jk}$ satisfy the following equation

$$\sigma_{jk} = \frac{a-cb}{(1-db)} \tag{120}$$

$$\omega_{jk} = \frac{c-ad}{1-bd} \tag{121}$$

Where,

$$a = \frac{1}{\gamma_1} - \frac{1}{g_{jk}(T) \, s'_{jk}(T) D' \, S^{-1}_{(y|u_{jk})}(T) D s_{jk}(T)} \tag{122}$$

$$b = \frac{d_{jk}(T) \, r'_{jk}(T) D' S^{-1}_{(y|u_{jk})}(T) \, D r_{jk}(T)}{g_{jk}(T) \, s'_{jk}(T) D' \, S^{-1}_{(y|u_{jk})}(T) D s_{jk}(T)} \tag{123}$$

$$c = \frac{1}{\gamma_2} - \frac{1}{d_{jk}(T) \, r'_{jk}(T) D' \, S^{-1}_{(y|u_{jk})}(T) D r_{jk}(T)} \tag{124}$$

$$d = \frac{g_{jk}(T) D s_{jk}(T) S^{-1}_{(y|u_{jk})} s'_{jk}(T) D'}{d_{jk}(T) \, r'_{jk}(T) D' \, S^{-1}_{(y|u_{jk})}(T) D r_{jk}(T)} \tag{125}$$

**proof:**

Since the matrices $G_1(b_k)$ and $G_2(b_{1k})$ are positive definite, then the matrix $G_T(b_k, b_{1k})$ is also positive definite for any k. Therefore, the controllability matrix $G_T$ is nonsingular. This implies that the system (4-7) is controllable from lemma (5.1).

Also since the decomposition of any positive definite matrix to a sum of one rank matrices is not unique therefore, if we take the eigenvalues decomposition of $G_1(b_k)$, we have

$$G_1(b_k) = S_k(T) O_k(T) S'_k(T) \tag{126}$$

Such that, $O_k(T)$ is a diagonal matrix of eigenvalues $g_{jk}(T)$, j=1,2,…,n of $G_1(b_k)$ and $S_k(T)S'_k(T) = I_{n \times n}$.

Then,
$$G_1(b_k) = \sum_{j=1}^{n} g_{jk}(T) \, s_{jk}(T) S'_k(T) \tag{127}$$

$s_{jk}(T)$, j=1,2,…,n are the columns of the matrix $S_k(T)$.

The matrix $G_1(b_k)$ is a positive definite by assumption, then, $g_{jk}(T) > 0$ for any j. If we choice $v_{jk}(T)$ as eigenvectors of the controllability matrix $G_1(b_k)$ then

$$v_{jk}(T) = \sqrt{g_{jk}(T)} s_{jk}(T), \text{ for every } j \in \{1,2,\dots,n\} \tag{128}$$

$$v_{jk}(T) = 0, \text{for every } j \in \{n+1, n+2, \dots\} \tag{129}$$

In similar manner, we take the eigenvalue decomposition of $G_2(b_{1k})$, then, we have

$$G_2(b_{1k}) = R_k(T) D_k(T) R'_k(T) \tag{130}$$

Where, $D_k(T)$ is a diagonal matrix of eigenvalues $d_{jk}(T)$, j=1,2,…,n of $G_2(b_{1k})$ and $R_k(T)R'_k(T) = I_{n \times n}$.

Then,
$$G_2(b_{1k}) = \sum_{j=1}^{n} d_{jk}(T) \, r_{jk}(T) r'_k(T) \tag{131}$$

$r_{jk}(T)$, j=1,2,…,n are the columns of a matrix $R_k(T)$.

But the matrix $G_2(b_{1k})$ is positive definite by assumption, then, $d_{jk}(T) > 0$ for any j. we choice $z_{jk}(T)$ as eigenvectors of the controllability matrix $G_2(b_{1k})$ then,

$$z_{jk}(T) = \sqrt{d_{jk}(T)} r_{jk}(T), \text{ for every } j \in \{1,2,\dots,n\} \tag{132}$$

$$z_{jk}(T) = 0, \text{ for every } j \in \{n+1, n+2, \dots\} \tag{133}$$

Taking the equations (128-129), (132-133) and equation (92) then expansion variances appear as:

For j=1,2,…,n and k=1,2,…,p

$$\sigma_{jk} = max\left\{0, \frac{1}{\gamma_1} - \frac{1}{g_{jk}(T) s'_{jk}(T) D' \, S^{-1}_{(y|u_{jk})}(T) D s_{jk}(T)} - \frac{\omega_{jk} d_{jk}(T) \, r'_{jk}(T) D' S^{-1}_{(y|u_{jk})}(T) \, D r_{jk}(T)}{g_{jk}(T) s'_{jk}(T) D' \, S^{-1}_{(y|u_{jk})}(T) D s_{jk}(T)}\right\} \tag{134}$$

and $\sigma_{jk}$ satisfies the power constraint (71). In addition, the equation (97) can be written as:

$$\omega_{jk} = max\left\{0, \frac{1}{\gamma_2} - \frac{1}{d_{jk}(T)\, r'_{jk}(T)D'\, S^{-1}_{(y|u_{jk})}(T)Dr_{jk}(T)} - \frac{\sigma_{jk}g_{jk}(T)Ds_{jk}(T)S^{-1}_{(y|u_{jk})}s'_{jk}(T)D'}{d_{jk}(T)\, r'_{jk}(T)D'\, S^{-1}_{(y|u_{jk})}(T)Dr_{jk}(T)}\right\} \quad (135)$$

$\omega_{jk}$ satisfies the constraint equation (72).

from assumption notations in this theorem the proof is completed.

### 7. Relationship between Capacity of Control and Final Time

In this section, we show that in general the capacity of control depends on the final time T. Precisely, we show that the capacity of control is bounded by the value of T in the following cases:

**First case**: infinite final time

The total covariance of uncontrolled noise matrix with the resulting matrix from initial random vector $x_0$, $S_{Total}(T)$ diverges for $T \to \infty$ the unlimited growth appears from the covariance matrix of noise $W^{H^2}$ with time T, which means that the capacity of control would go to zero as $T \to \infty$. We found that the capacity of control is finite both in asymptotically stable systems and in unstable system (all the eigenvalues of A are positive real part).

i. **Stable systems:**

The control system (4-7) is stable, when the total covariance of uncontrolled noise matrix $S_{Total}(T)$ and controllability matrix $G_T$ are finite for $T \to \infty$.

By solving the following Lyapunov equations for any k = 1, 2, …,p

$$AG_1(b_k, \infty) + G_1(b_k, \infty)A' = -b_k b_k' \quad (136)$$

$$AS_{dW}(\infty) + S_{dW}(\infty)A' = G\ G' \quad (137)$$

$$AS_{dW^{H^1}} + S_{dW^{H^1}}A' = \sigma_1\ \sigma_1' \quad (138)$$

We get,

$$\sigma_{jk} = max\left\{0, \frac{1}{\gamma_1} - \frac{1}{g_{jk}(\infty)\, s'_{jk}(\infty)D'\, S^{-1}_{(y|u_{jk})}(\infty)Ds_{jk}(\infty)} - \frac{\omega_{jk}d_{jk}(\infty)\, r'_{jk}(\infty)D'S^{-1}_{(y|u_{jk})}(\infty)\ Dr_{jk}(\infty)}{g_{jk}(\infty)\, s'_{jk}(\infty)D'\, S^{-1}_{(y|u_{jk})}(\infty)Ds_{jk}(\infty)}\right\} \quad (139)$$

$$\omega_{jk} = max\left\{0, \frac{1}{\gamma_2} - \frac{1}{d_{jk}(\infty)\, r'_{jk}(\infty)D'\, S^{-1}_{(y|u_{jk})}(\infty)Dr_{jk}(\infty)} - \frac{\sigma_{jk}g_{jk}(\infty)Ds_{jk}(\infty)S^{-1}_{(y|u_{jk})}s'_{jk}(\infty)D'}{d_{jk}(\infty)\, r'_{jk}(\infty)D'\, S^{-1}_{(y|u_{jk})}(\infty)Dr_{jk}(\infty)}\right\} \quad (140)$$

This implies that, the capacity of control $C(\infty)$ is finite for any power constrains $M_1$ and $M_2$, where the total covariance matrix and the controllability matrix $G_T$ completely defines the capacity of control.

2) **Unstable systems**

In this case, we show that the capacity of control is finite as well, if all eigenvalues of a matrix A are positive real part.

**Lemma (7.1)**

The following controllability matrices and process noise variance matrix respectively,

$$\widetilde{G}_1(b_{1k}) = \int_0^{T-h} e^{-At} b_k b_k' e^{-A't}\, dt \quad (141)$$

$$\widetilde{G}_2(b_{1k}) = \int_{T-h}^{T} e^{-At} b_{1k} b_{1k}' e^{-A't}\, dt \quad (142)$$

$$\widetilde{S}_{Total}(T) = DC_{x_0}D' + D\int_0^T e^{-At} G\ G'\ e^{-A't}\, dt\, D' + D\, H\,(2H-1)\int_0^T e^{-At} \sigma_1 |T-t|^{2H-2}\, \sigma_1' e^{-A't}\, dt\ D'$$

$$+\sigma_2 \begin{bmatrix} T^{2\widetilde{h}_1} & & 0 \\ & \ddots & \\ 0 & & T^{2\widetilde{h}_P} \end{bmatrix}_{P\times P} \sigma'_2 \quad (143)$$

Define the same mutual information objective function in (60), as the original matrices $G_1(b_k)\ G_2(b_{1k})$ and $S_{Total}(T)$.

**Proof:**

Form the optimal sets of $\{v_{jk}(T)\}$ and $\{z_{jk}(T)\}$ the objective function (60) can be written as

$$I(y(T); u(t), t \in [0, T]) = \ln\{\det[I_{n\times n} + S_{Total}^{-1}(T) \sum_{k=1}^{p} D\, G_1(b_K)D' + S_{Total}^{-1}(T) \sum_{k=1}^{p} D\, G_2(b_{1^k})D']\} \quad (144)$$

where, $S_{Total}^{-1}(T)$, $G_1(b_K)$ and $G_2(b_{1^k})$ can be represented by:

$$S_{Total}^{-1}(T) = [e^{A'T}]^{-1} (\tilde{S}_{Total}(T))^{-1} [e^{AT}]^{-1}$$

$$G_1(b_K) = [e^{AT}]\tilde{G}_1(b_K)[e^{A'T}]$$

$$G_2(b_{1^k}) = [e^{AT}]\tilde{G}_2(b_{1^k})[e^{A'T}]$$

Therefore, the expression in (144) is equivalent to:

$$I(y(T); u(t), t \in [0, T]) = \ln\left\{\det\left[I_{n\times n} + (\tilde{S}_{Total}(T))^{-1} \sum_{k=1}^{p} D\tilde{G}_1(b_k)D' + (\tilde{S}_{Total}(T))^{-1} \sum_{k=1}^{p} D\tilde{G}_2(b_{1^k})D'\right]\right\} \quad (145)$$

When all eigenvalues of A are positive (the system is unstable) then, the following equations are hold

$$A\tilde{S}_{dW}(\infty) + \tilde{S}_{dW}(\infty)A' = G\, G' \quad (146)$$

$$A\tilde{S}_{dW^{H^1}} + \tilde{S}_{dW^{H^1}}A' = \sigma_1\, \sigma_1' \quad (147)$$

$$A\tilde{G}_1(b_k, \infty) + \tilde{G}_1(b_k, \infty)A' = -b_k b_k' \quad (148)$$

$$A\tilde{G}_2(b_{1^k}, \infty) + \tilde{G}_2(b_{1^k}, \infty)A' = -b_{1^k} b_{1^k}' \quad (149)$$

By solving the corresponding Lyapunov equations above to get the solution by replacing

$$G_1(b_k, \infty) = \tilde{G}_1(b_k, \infty)$$

$$S_{Total}(T) = \tilde{S}_{Total}(T)$$

**Second case:** infinitesimal final time

In this case, we show that the capacity of control becomes a linear function of final time T, when T goes to zero with constant coefficient that depends on its parameters of the system. In particular the capacity of control is vanished for $T \to 0$. By approximating integral in equation (24) for $T \to 0$, we get

$$S_{Total}(T) = (DG\, G'D' + D\sigma_1\sigma_1'D' + \sigma_2\sigma_2')T \quad (150)$$

Similarity, for T takes values approaching to zero, the effect of h is not almost existent, so it can be assumed zero. Then the control process covariance matrix $S_u(T)$ can be written as:

$$S_u(T) = \sum_{k=1}^{p} \sum_{j=1}^{n} \sigma_{jk}\, b_k e_{j,k}(0) b_k' e_{j,k}(0)\, T^2 \quad (151)$$

Therefore, the capacity of control for T near than zero appears as

$$C(T) = \ln\{\det[I_{n\times n} + S_{Total}^{-1}(T) \sum_{k=1}^{p} D\, G_1(b_K)D' + S_{Total}^{-1}(T) \sum_{k=1}^{p} D\, G_2(b_{1^k})D']\} \quad (152)$$

Substituting (150-151) in equation (152), we get

$$C(T) = \ln\{\det[I_{n\times n} + Q_{n\times n}T]\} \quad (153)$$

Consequently,

$$C(T) = \text{Tr}\{Q_{n\times n}\}\, T \quad (154)$$

where, $Q_{n\times n} = (DG\, G'D' + D\sigma_1\sigma_1'D' + \sigma_2\sigma_2')^{-1} \cdot \sum_{k=1}^{p} \sum_{j=1}^{n} \sigma_{jk}\, b_k e_{j,k}(0) b_k' e_{j,k}(0)$ is constant matrix.

### 7. Illustrative example

As a simple illustrative example, explains the relationship between capacity of control and the intrinsic properties of stochastic dynamic systems perturbed by mixed fractional Brownian motion with delay in control, consider control system of the form (4-7) defined in a given time interval [0,T], assume that $T > \mathbf{h}$.

**Example (7.1):**

Consider the following constant matrices

$$A = \begin{bmatrix} -1 & 0 & 1 \\ 0 & -2 & 0 \\ 0 & 0 & -1 \end{bmatrix}, B_1 = \begin{bmatrix} 0 & 1 \\ 1 & 1 \\ 1 & -1 \end{bmatrix}, B_2 = \begin{bmatrix} 0 & 0 \\ 1 & 0 \\ 0 & 1 \end{bmatrix}, G = \begin{bmatrix} 1 & 0 & 0 \\ 0 & 1 & 0 \\ 0 & 0 & 1 \end{bmatrix}, \sigma_1 = \begin{bmatrix} 1 & 1 & 0 \\ -1 & 1 & 0 \\ 2 & 0 & 1 \end{bmatrix},$$

$$D = \begin{bmatrix} 1 & -1 & 0 \\ -3 & 1 & 0 \\ 2 & 0 & 1 \end{bmatrix}, \sigma_2 = \begin{bmatrix} 1 & 0 & 0 \\ 0 & 1 & 0 \\ 0 & 0 & 1 \end{bmatrix}, H_1 = \begin{bmatrix} 0.75 \\ 0.75 \\ 0.75 \end{bmatrix}, H_2 = \begin{bmatrix} 0.75 \\ 0.75 \\ 0.75 \end{bmatrix}, Cx0 = \begin{bmatrix} 1 & 0 & 0 \\ 0 & 2 & 0 \\ 0 & 0 & 1 \end{bmatrix}$$

Hence, n=3, p=2 and

assume that T= 2 with constant point delay h=1. Moreover, $B = B_1 + e^{-Ah}B_2 = \begin{bmatrix} 1 & 1.3679 \\ 8.389056 & 2.0000 \\ 2 & 1.7183 \end{bmatrix}$

If $M_1 + M_2 \leq 1$ then $G_1(b_1) = \begin{bmatrix} 329.4882 & 62.0411 & 176.6757 \\ 62.0411 & 12.2614 & 34.4356 \\ 176.6757 & 34.4356 & 97.1219 \end{bmatrix}$,

$$G_1(b_2) = \begin{bmatrix} 111.8777 & 32.9990 & 38.7258 \\ 32.9990 & 10.2647 & 12.0753 \\ 38.7258 & 12.0753 & 14.2068 \end{bmatrix}, G_2(b_{1^1}) = \begin{bmatrix} 7.6311 & 3.7827 & 4.3504 \\ 3.7827 & 2.1101 & 2.3812 \\ 4.3504 & 2.3812 & 2.6966 \end{bmatrix}$$

$$G_2(b_{1^2}) = \begin{bmatrix} 0.4545 & 0.1169 & -0.2368 \\ 0.1169 & 0.2454 & 0.5479 \\ -0.2368 & 0.5479 & 1.9039 \end{bmatrix}$$

Note that the matrices $G_1(b_1)$, $G_1(b_2)$, $G_2(b_{1^1})$ and $G_2(b_{1^2})$ are positive definite, then the controllability matrix $G_T$ is nonsingular. Hence, by lemma (5.1), the system is controllable, this mean that the capacity of control not equal zero. The optimal vectors are

$$v_{11}(T) == \begin{bmatrix} 18.1368 \\ 3.4458 \\ 9.7879 \end{bmatrix} \quad v_{21}(T) = \begin{bmatrix} 0.7370 \\ -0.6188 \\ -1.1478 \end{bmatrix}, v_{31}(T) = \begin{bmatrix} 0.0051 \\ 0.0683 \\ -0.0335 \end{bmatrix} \quad v_{12}(T) == \begin{bmatrix} 10.5681 \\ 3.1473 \\ 3.6952 \end{bmatrix}$$

$$v_{22}(T) = \begin{bmatrix} 0.4383 \\ -0.5991 \\ -0.7431 \end{bmatrix} \quad v_{32}(T) = \begin{bmatrix} 0.0000122 \\ -0.0009262 \\ 0.0007539 \end{bmatrix}, z_{11}(T) == \begin{bmatrix} 2.7430 \\ 1.4170 \\ 1.6188 \end{bmatrix}, z_{21}(T) = \begin{bmatrix} 0.3271 \\ -0.3189 \\ -0.2750 \end{bmatrix}$$

$$z_{31}(T) = \begin{bmatrix} 0.0020 \\ 0.0200 \\ -0.0208 \end{bmatrix} \quad z_{12}(T) == \begin{bmatrix} 0.1711 \\ -0.3986 \\ -1.3796 \end{bmatrix}, z_{22}(T) = \begin{bmatrix} 0.6515 \\ 0.2870 \\ -0.0021 \end{bmatrix}, z_{32}(T) = \begin{bmatrix} 0.0285 \\ -0.0646 \\ 0.0222 \end{bmatrix}$$

$$S_{Total}(T) = \begin{bmatrix} 77.6392 & -235.9072 & 170.4309 \\ -235.9072 & 825.5392 & -619.2593 \\ 170.4309 & -619.2593 & -619.2593 \end{bmatrix}$$

Capacity = 1.9071

Example (7.2):

Consider the control system of the form (1) is defined in a given time interval [0,T]. Assume that T > **0** with the following constant matrices

$A=\begin{bmatrix}-1 & 0\\ 0 & -2\end{bmatrix}$ , $B_1=\begin{bmatrix}1\\0\end{bmatrix}$ , $B_2=\begin{bmatrix}0\\1\end{bmatrix}$ , $G=\begin{bmatrix}1 & 0\\ 0 & 1\end{bmatrix}$ , $\sigma_1=\begin{bmatrix}1 & 0\\ 1 & -1\end{bmatrix}$ , $D=\begin{bmatrix}3 & 0\\ 0 & 1\end{bmatrix}$ , $\sigma_2=\begin{bmatrix}1 & 0\\ 0 & 1\end{bmatrix}$

$H_1=\begin{bmatrix}0.75\\0.75\end{bmatrix}$ , $H_2=\begin{bmatrix}0.75\\0.75\end{bmatrix}$ , $Cx0=\begin{bmatrix}1 & 0\\ 0 & 2\end{bmatrix}$

Hence, n=2, p=1

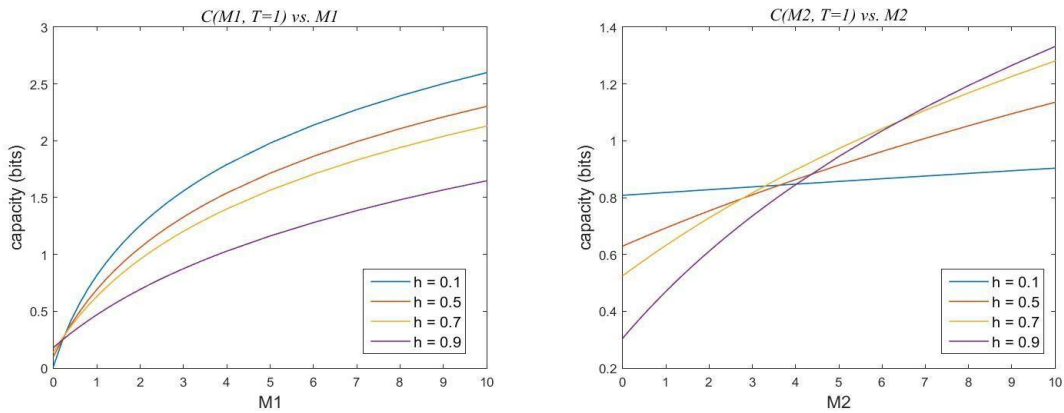

Fig. 1 The horizontal line represents the power constraints M1, M2 respectively and the vertical line is the entropy of observe state at T=1.

As seen in figure 1, the graphs of capacity of control $C_h(M_1, M_2=1, T=1)$ become more similar to each other, when the time delay h decreases. In contrast working when the time delay h increases, the capacity of control $C_h(M_2, M_1=1, T=1)$ become more similar to each other. Also, we observe that the capacity of control $C_h(M_2, M_1=1, T=1)$ tends approximate to be constant for any $M_2$ when h converges to zero.

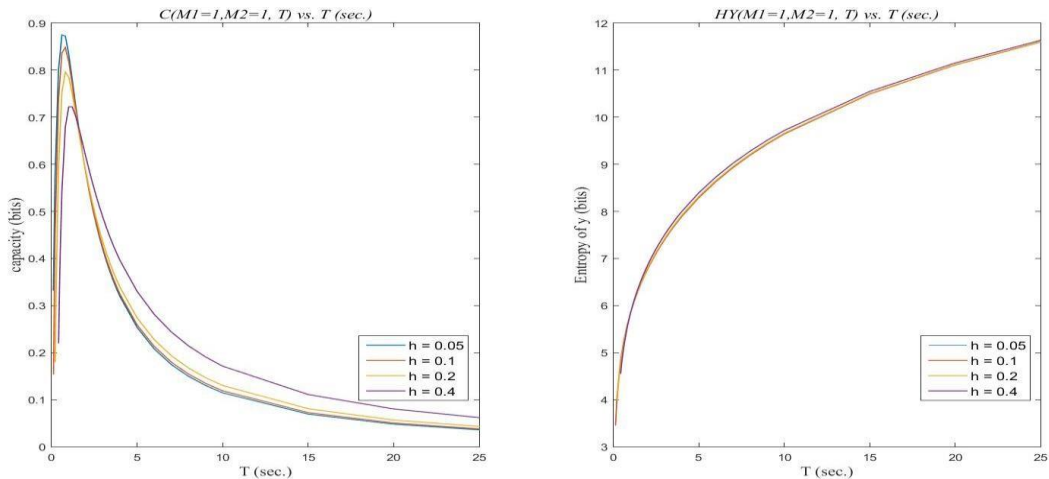

Fig. 2 the horizontal line represents the time and the vertical lines are the capacity of control and the entropy of observe state, respectively.

As shown in Fig. 2, for any T, the capacity of control $C_h(M_1=1, M_2=1, T)$ converges to the finite value, for decreasing h. This behavior is qualitatively similar for different power constraints $M_1$ and $M_2$ of dynamical system. Also, we observe that, there is a relationship between the value of capacity of control and the values of T and h, $C_h(M_1=1, M_2=1, T)$ converges to 0.9118 for h converges to zero this is true when T=0.6300 on the other hand, as always the entropy of observe state $H(y(T))$ is increasing when T increasing. Mathematically, the ratio of the noise matrix to the controllability matrix in (145) converges to a constant for T converges to zero.

.

.